\documentclass[12pt]{article}
\usepackage{amsmath, amsthm, amssymb, amsfonts, enumerate}
\usepackage[figuresright]{rotating}
\usepackage{natbib}

\begin{document}
\newcommand{\bea}{\begin{eqnarray}}
\newcommand{\ena}{\end{eqnarray}}
\newcommand{\beas}{\begin{eqnarray*}}
\newcommand{\enas}{\end{eqnarray*}}
\newcommand{\beq}{\begin{equation}}
\newcommand{\enq}{\end{equation}}
\def\qed{\hfill \mbox{\rule{0.5em}{0.5em}}}
\newcommand{\bbox}{\hfill $\Box$}
\newcommand{\From}{From}
\newcommand{\ignore}[1]{}
\newcommand{\ws}{{\widetilde \sigma}}
\newcommand{\btheta}{\mbox{\boldmath {$\theta$}}}
\newcommand{\transpose}{{\mbox{\scriptsize\sf T}}}
\newcommand{\esssup}{\mathop{\mathrm{ess\,sup}}}
\newcommand{\B}{B}
\newcommand{\D}{D}
\newcommand{\E}{E}
\newcommand{\V}{V}
\newcommand{\Km}{K_{\eta,m}}
\newcommand{\Eta}{\eta}
\newcommand{\de}{\varrho}
\newcommand{\tP}{\overline{P}}
\newcommand{\tD}{\overline{D}}
\newcommand{\tQ}{\overline{Q}}
\newcommand{\tL}{\overline{L}}
\newcommand{\tN}{\overline{N}}
\newcommand{\fmx}{\sup}
\newtheorem{theorem}{Theorem}[section]
\newtheorem{corollary}{Corollary}[section]
\newtheorem{conjecture}{Conjecture}[section]
\newtheorem{proposition}{Proposition}[section]
\newtheorem{lemma}{Lemma}[section]
\newtheorem{definition}{Definition}[section]
\newtheorem{example}{Example}[section]
\newtheorem{remark}{Remark}[section]
\newtheorem{case}{Case}[section]
\newtheorem{condition}{Condition}[section]
\newtheorem{assumption}{Asssumption}[section]

\newcommand{\To}{\rightarrow}
\newcommand{\wtilde}[1]{\widetilde{#1}}
\newcommand{\qmq}[1]{\quad\mbox{#1}\quad}
\newcommand{\changed}{\texttt{[CHANGED:]}}

\newcommand{\figformat}{pdf}

\title{{\bf\Large On Optimal Allocation of a Continuous Resource Using an Iterative Approach
and Total Positivity}}

\author{Jay Bartroff$^1$, Larry Goldstein$^1$, Yosef Rinott$^2$ and Ester Samuel-Cahn$^3$\\\\
\small{$^1$ Department of Mathematics, University of Southern California, Los Angeles, California, USA}\\
\small{$^2$ Department of Statistics and Center for the Study of
Rationality,}\\ \small{The Hebrew University of Jerusalem, Israel, and LUISS, Rome}\\
\small{$^3$ Department of Statistics and Center for the Study of
Rationality,}\\ \small{The Hebrew University of Jerusalem,
Israel}\\} \footnotetext{AMS 2000 subject classifications. Primary
60G40\ignore{Stopping times; optimal stopping problems; gambling
theory}, Secondary 62L05\ignore{Sequential design},
91A60\ignore{optimal stopping times, sequential decision,
Probabilistic games }} \footnotetext{Key words and phrases: bomber
problem,  log-concave, optimal allocation, sequential optimization,
total positivity} \maketitle

\date{}

\begin{abstract}
We study a class of optimal allocation problems, including the
well-known Bomber Problem, with the following common probabilistic
structure. An aircraft equipped with an amount~$x$ of ammunition is
intercepted by enemy airplanes arriving according to a homogenous
Poisson process over a fixed time duration~$t$. Upon encountering an
enemy, the aircraft has the choice of spending any amount~$0\le y\le x$
of its ammunition, resulting in the aircraft's survival with probability equal to some known increasing function of $y$.
Two different goals have been considered in the literature
concerning the optimal amount~$K(x,t)$ of ammunition spent: (i)~Maximizing the probability of surviving for time~$t$, which is the
so-called Bomber Problem, and (ii) maximizing the number of enemy
airplanes shot down during time~$t$, which we call the Fighter
Problem. Several authors have attempted to settle the following
conjectures about the monotonicity of $K(x,t)$: [A] $K(x,t)$ is
decreasing in $t$, [B] $K(x,t)$ is increasing in $x$, and
[C] the amount~$x-K(x,t)$ held back is increasing in $x$. [A]
and [C] have been shown for the Bomber Problem with discrete
ammunition, while [B] is still an open question. In this paper we consider both time and ammunition continuous, and
for the Bomber Problem prove [A] and [C], while for the Fighter
we prove [A] and [C] for one special case and
[B] and [C] for another. These proofs  involve showing that the optimal survival probability
and optimal number shot down are totally positive of order 2
($\mbox{TP}_2$) in the Bomber and Fighter Problems, respectively.
The $\mbox{TP}_2$ property is shown by constructing convergent
sequences of approximating functions through an iterative operation
which preserves $\mbox{TP}_2$ and other properties.\end{abstract}

\section{Introduction}


Throughout this paper we write decreasing for non-increasing and increasing for
non-decreasing. When we use strict monotonicity, we will say so explicitly.
The allocation problem  discussed in this paper can be described in
terms of the following example: an aircraft, which is at flying time~$t$ away from its final destination, is equipped with an amount $x$
of ammunition. It is confronted by enemy airplanes whose appearance
is driven by a time-homogenous Poisson process with known intensity.
By adjusting the units of time, we assume without loss of generality
that the intensity is one. It is assumed that the encounters are
instantaneous, that is, that their duration time is zero. The
question of interest is how much of the presently available
ammunition should be spent when confronted by an enemy while in
``state''~$(x,t)$. We consider two problems, the Bomber Problem and
the Fighter Problem, with this common probabilistic setup. For the
Bomber Problem, the goal is to maximize the probability~$P(x,t)$ of
reaching its destination. For the Fighter Problem, the goal  is to
maximize the expected number~$N(x,t)$ of enemy airplanes shot down. In both problems there is given an increasing function~$a(y)$
which, in the Bomber Problem is the probability that the Bomber
\textit{survives} an enemy encounter in which it expends $0\le y\le x$
units of its ammunition, while in the Fighter Problem $a(y)$ is the
probability that the Fighter \textit{destroys} an enemy when it
expends $y$ units of its ammunition. Although our results below are
proved when $a(y)$ is general, the canonical example considered in the
literature for $a(y)$ in the Bomber Problem is
\begin{equation}\label{a-function}
1-(1-u)e^{-y}
\end{equation}
for some fixed $u\in[0,1]$, which can be interpreted as the Bomber's
$y$ units of ammunition destroying the enemy with
probability~$1-e^{-y}$, while otherwise allowing the enemy to launch a counterattack
which succeeds with probability~$1-u$.

For the Fighter, where $a(y)$ is the probability of destroying the enemy,
the probability that the Fighter survives an encounter is $a(y)
+u(1-a(y))$, since the enemy is hit with probability $a(y)$, and
otherwise, with probability $1-a(y)$ the enemy counterattacks, and misses
with probability $u$. When taking the canonical choice $a(y)=1-e^{-y}$,
the Fighter's survival is given by  (\ref{a-function}).

While the optimal ammunition spending strategies
$K(x,t)$ may differ for the Bomber and Fighter Problems, it is
perhaps intuitively obvious that the following three monotonicity
properties, first posed by \citet{Klinger68} for the Bomber Problem,
should hold for both:
\begin{enumerate}

\item[][A]: $K(x,t)$ is decreasing in $t$.

\item[][B]: $K(x,t)$ is increasing in $x$.

\item[][C]: $x-K(x,t)$ is increasing in $x$.

\end{enumerate}
Roughly speaking, [A] states that the closer one is to the
destination, the more one should spend, [B] states that the more
ammunition one has, the more one should spend, and [C] states that
the more ammunition one has, the more one should retain for possible
future encounters; replacing ``more'' by ``not less'' makes these
characterizations precise. Although we write $K(x,t)$ for the two
different goals, no confusion should arise, as it will always be
clear which of the problems we are discussing, or that we are
discussing them simultaneously. An optimal policy may not be unique.
If it is not unique, $K(x,t)$ will always denote the minimal optimal
amount of ammunition needed.

Although the paper is written in \textit{plane} English, the
settings above serve only as illustrations of the general method which can be applied to a
wide class of problems where some limited resource needs to be optimally allocated over time or space to maximize the probability of a system's survival, which clearly also have more peaceful applications.

Of the two problems, the Bomber Problem has received most attention. For discrete ammunition, such as missiles, it is first treated in \citet{Klinger68}, who prove that [B] implies [A]. A clearer picture emerges in \citet{Samuel70}, who proves [A] without assuming [B]. Samuel's \citeyearpar{Samuel70} formula (3.9) also proves [C]. Conjecture [B] has remained elusive, and in fact remains in doubt. \citet{Simons90} argue that when both time and ammunition are continuous,  then [B] is equivalent
to the log concavity of $P(x,t)$ in $x$ for each fixed $t$. They also show that when both time
and ammunition are discrete, there exist parameter
values for which log concavity fails. However, they are unable to supply an
example where [B] fails.

Regarding the Fighter Problem, we are unaware
of this variation on the Bomber Problem being discussed in a finite time~$t$ horizon setting.
For an infinite time horizon, an optimal policy can be written simply as $K(x)$.  \citet{Weber85} shows that the monotonicity of $K(x)$ fails to hold for discrete ammunition when $a$ is given by (\ref{a-function}) with $u=0$, and \citet{Shepp91} show that monotonicty fails for continuous ammunition.  In a setting where both ammunition and time are discrete, the duration $T$ is an unknown exponential random variable with known mean, and geometric arrivals of the enemy, \citet{Shepp91} also show that there exist parameter values for which [B] fails. Though this variant is on a finite time interval of length $T$, the optimal policy does not depend on $t$ and so resembles an infinite time horizon problem. On the other hand, \citet{Bartroff10} show for the Bomber Problem with continuous time and ammunition that there exists a certain  region of the $(x,t)$-space over which $K(x,t)=x$,  a monotone function of $x$, and hence  [B] holds in this region.

In the present paper we consider the ``doubly continuous'' setting
where ammunition and time are continuous variables. Given an
increasing function~$a(y)$ taking values in $[0,1]$, the dynamic
programming equation for the Bomber Problem giving the optimal
survival probability~$P(x,t)$ is
$$
P(x,t)=\int_0^t  \sup_{0 \le y \le
x}a(y)P(x-y,t-s)e^{-s} ds+e^{-t} \quad \mbox{for $(x,t) \in
\mathbb{R}^+ \times \mathbb{R}^+$},
$$ where $\mathbb{R}^+=[0,\infty).$ The explanation is as follows. There are no encounters in a time
interval of length~$t$ with probability~$e^{-t}$, yielding the final
term. Otherwise, an encounter occurs at time $s \in [0,t]$ with
density~$e^{-s}$, at which time the amount $y \in [0,x]$ will be
chosen to maximize the probability~$a(y)$ of surviving the current
encounter times the probability of future survival under an
optimal policy. A simple change of variables yields
\begin{equation}\label{P}
P(x,t)=\int_0^t \sup_{0 \le y \le
x} a(y)P(x-y,s)e^{-(t-s)} ds+e^{-t} \quad \mbox{for $(x,t) \in
\mathbb{R}^+ \times \mathbb{R}^+$}.
\end{equation}

Given an increasing function~$a(y)$ taking values in $[0,1]$ and
a fixed value~$u\in[0,1]$, the dynamic programming equation for the
Fighter Problem giving the optimal expected number~$N(x,t)$ of
enemies shot down is
\begin{equation}\label{eq:gLPI}
N(x,t)=\int_0^t
\sup_{0 \le y \le x}\{a(y)+[a(y) +u(1-a(y))]N(x-y,s)\}e^{-(t-s)}ds \quad
\mbox{for $(x,t) \in \mathbb{R}^+ \times \mathbb{R}^+$},
\end{equation}
which can be interpreted as follows. Encounters occur with density
$e^{-s}$, and expending amount~$y$ the Fighter gains a single hit if
he destroys the enemy, which happens with probability~$a(y)$. The
Fighter will gain  an expected additional $N(x-y,t-s)$ future
hits if he destroys the enemy or otherwise if the enemy's
counterattack fails, which happens with probability~$u$. A change of
variables then yields (\ref{eq:gLPI}).

By putting both the Bomber and the Fighter in a common general framework,
we show that the solutions to (\ref{P}) and (\ref{eq:gLPI})
exist and are unique and continuous. In particular, the supremum is
attained in both (\ref{P}) and (\ref{eq:gLPI}), and hence
$\sup$ can be replaced by $\max$ in each.
The optimal policy~$K(x,t)$ is then defined as the minimal value for
which the $\max$ in the integrands in (\ref{P}) and (\ref{eq:gLPI})
is attained, respectively.

Recall that a nonnegative function $Q(x,t)$ is \textit{totally positive of order~2}, written $\mbox{TP}_2(x,t)$ or simply $\mbox{TP}_2$, if
\begin{equation}
\label{TP2-def}
Q(x',t')Q(x,t) \ge Q(x',t)Q(x,t') \quad \mbox{whenever $x<x'$ and $t<t'$,} \end{equation}
and strictly $\mbox{TP}_2$ if $\ge$ in (\ref{TP2-def}) can be replaced by $>$. (See \citet{Karlin68} and the references therein.) To prove [A] for the Bomber Problem, we establish the TP$_2$ property of $P(x,t)$ in Section~\ref{sec:bomber}. Recall also that a \textit{log-concave function} is a function whose logarithm is concave.

\begin{theorem}
\label{main} Let $a(y)$ be a uniformly continuous log-concave function. Then a unique bounded solution $P(x,t)$ to equation (\ref{P})
exists, satisfies $P(x,t) \in [0,1]$, is continuous and $\mbox{TP}_2$.
\end{theorem}

\begin{theorem}
\label{thm:A} If $a(y)$ is a uniformly continuous log-concave function, then [A] holds for the doubly continuous Bomber Problem.
\end{theorem}

\noindent\citet{Simons90} claim that [C] can be shown to hold for continuous $x$ and $t$ by arguments similar to the ones they provide for a case where both $x$ and $t$ are discrete; in Section \ref{sec:bomber} we provide a rigorous proof of [C] for the doubly continuous case.
The discrete $x$ analogues to
Theorems~\ref{main} and \ref{thm:A} were proved in \citet{Samuel70}.

Regarding the Fighter in general, we have the following result.
\begin{theorem}
\label{thm:genFighter}
If $a(y)$ is uniformly continuous, then for any $u \in [0,1]$ a bounded solution $N(x,t)$ to equation (\ref{eq:gLPI}) exists, is unique, satisfies $N(x,t) \in [0,t]$, and is continuous.
\end{theorem}

For the Fighter Problem we do not resolve the question of whether the general solution to (\ref{eq:gLPI}) is $\mbox{TP}_2$, but instead examine two special cases in more detail:

i) \textit{The Frail Fighter}: the case $u=0$. Once the Frail Fighter fails to shoot down an enemy, he himself is shot down with probability one. Equation~(\ref{eq:gLPI}) then simplifies to
\begin{equation}
\label{eq:LPI} N(x,t)=\int_0^t \sup_{0 \le
y \le x}a(y)[1+N(x-y,s)]e^{-(t-s)}ds \quad \mbox{for $(x,t) \in
\mathbb{R}^+ \times \mathbb{R}^+$}.
\end{equation}

\begin{theorem}
\label{thm:1.3}
If $a(y)$ is a uniformly continuous log-concave function, then the solution $N(x,t)$ is $\mbox{TP}_2$, and property [A] holds for the Frail Fighter.
\end{theorem}

\noindent The Frail Fighter is discussed in Section \ref{sec:Ffighter}, where [C] is also established.

ii) \textit{The Invincible Fighter}: the case $u=1$. As the name suggests, the Invincible Fighter cannot be shot down. Equation (\ref{eq:gLPI}) simplifies to
\begin{equation}
\label{eq:DPI} N(x,t)=\int_0^t \sup_{0 \le y \le
x}[a(y)+N(x-y,s)]e^{-(t-s)}ds \quad \mbox{for $(x,t) \in
\mathbb{R}^+ \times \mathbb{R}^+$}.
\end{equation}
\begin{theorem}
\label{thm:1.4}
If $a(y)$ is a uniformly continuous concave function then
property [B] holds for the Invincible Fighter.
\end{theorem}

\noindent The Invincible Fighter is discussed in Section
\ref{sec:Ifighter}, where we also show [C] and that the maximum $y$
determining $K(x,t)$ is uniquely attained. It should be noted that
for none of the problems considered here have we managed to show
that [A] and [B] hold simultaneously.

While in the discrete ammunition case one can demonstrate claims by
using induction on the number of available units \citep[e.g.,][]
{Simons90}, here we must take a different approach. In Section
\ref{sec:contraction} we prove the existence and uniqueness of
solutions to equations which have the form of (\ref{P}) or
(\ref{eq:gLPI}) by following the outline suggested by Weber
\citep[p.~431]{Simons90} and construct a sequence of functions
through an iteration which converges exponentially fast to the
solution.

We illustrate our approach using the Bomber Problem;
similar remarks apply to the Fighter. First, consider a
slightly transformed version of (\ref{P}). For ease of notation, for
functions $a(\cdot)$ and $Q(\cdot,\cdot)$ let
\begin{equation}\label{eq:star}
(a \star Q)(x,s)= \sup_{0 \le y \le x}a(y)Q(x-y,s),
\end{equation}
the sup-convolution of $a(\cdot)$ and $Q(\cdot,s)$. Letting
\begin{equation}
\label{tP=etP} \tP(x,t)=e^tP(x,t),
 \end{equation} substitution into (\ref{P})
results in the integral equation \begin{equation}
\label{tP-integral-equation} \tP(x,t)=\int_0^t (a \star
\tP)(x,s)ds+1 \quad \mbox{for $(x,t) \in \mathbb{R}^+ \times
\mathbb{R}^+$.} \end{equation} Clearly  there exists a unique,
continuous solution $\tP$ to (\ref{tP-integral-equation}) satisfying
$0 \le \tP(x,t) \le e^t$ if and only if there exists a unique,
continuous solution $P$ to (\ref{P}) satisfying $0 \le P(x,t) \le
1$; Theorem \ref{thm:Pcmt} proves the former and provides an
iteration method to approximate the solution. Theorem
\ref{thm:main2} provides further properties of the solution.

To prove the claims of Theorem \ref{thm:main2}, consider the
sequence of functions $\{\tP_m\}_{m\ge 0}$ generated with some
initial $\tP_0(x,t)$ by the recursion \begin{equation}
\label{recursion} \tP_{m+1}(x,t)=\int_0^t (a \star \tP_m)(x,s)ds+1
\quad \mbox{for} \quad m \ge 0. \end{equation} In Section
\ref{sec:contraction}, by showing that the operation which gives
$\tP_{m+1}$ from $\tP_m$ is a type of contraction, we are able to
prove the existence and uniqueness of the solution to
(\ref{recursion}) by the method of the contraction mapping theorem.
In Section \ref{sec:bomber} we prove that various operations
preserve the $\mbox{TP}_2$ and other properties, allowing us to
prove inductively that $\tP_m$ possesses a particular property for
all~$m\ge 0$, that may then be inherited by $\tP$. We remark in
Section \ref{sec:disc} why this same approach appears to fail to
give [B]. Nevertheless, the sequence of approximate solutions, which
are guaranteed to converge in the supremum norm at an exponential
rate by Theorem \ref{thm:Pcmt}, may yet have uses in addition to
that of proving [A].

\section{Existence, uniqueness, and approximating sequences}
\label{sec:contraction} We will put our models into the following
common framework. Given a function $a:\mathbb{R}^+
\rightarrow \mathbb{R}^+$, we consider an operation $\otimes$ that, for a given  function $Q:\mathbb{R}^+ \times \mathbb{R}^+
\rightarrow \mathbb{R}^+$,  returns a function $a
\otimes Q:\mathbb{R}^+ \times \mathbb{R}^+ \rightarrow
\mathbb{R}^+$.  Although our results in this section are more
general, we will be interested in the specific cases where
$(a\otimes Q)(x,t)$ is one of the following: \bea \label{twostars}
(Bo)&& \sup_{0 \le y \le x}a(y)Q(x-y,t),\nonumber\\
(F_0)&&\sup_{0 \le y \le x}a(y)[e^t
+Q(x-y,t)] \nonumber\\
(F_1)&&\sup_{0 \le y \le x}[ a(y)e^t+Q(x-y,t)], \nonumber \\
(F_u)&&\sup_{0 \le y \le
x}\{a(y)e^{t}+[a(y)+u(1-a(y))]{Q}(x-y,t)\}. \ena Case~$(Bo)$ corresponds
to the Bomber Problem \eqref{tP-integral-equation}, case~$(F_u)$ to the general Fighter
equation \eqref{eq:gLPI} after a rescaling by $e^t$ as in (\ref{tP=etP}), ~$(F_0)$ the case $u=0$ of ~$(F_u)$, the Frail Fighter, and case~$(F_1)$ the case $u=1$ of $(F_u)$, the Invincible Fighter.

Given $a$ and $\otimes$ we will consider operators mapping the collection of functions $Q:\mathbb{R}^+ \times
\mathbb{R}^+ \rightarrow \mathbb{R}^+$ to itself. For the Fighter Problem, these operators will be of the form
\begin{equation}
\label{def:calG}
{\cal G}(x,t,Q)=\int_0^t (a \otimes Q)(x,s)ds,
\end{equation}
and for the Bomber Problem,
\begin{equation}
\label{def:calFG} {\cal F}(x,t,Q)=\int_0^t (a \otimes Q)(x,s)ds \qmq{and} {\cal G}(x,t,Q)={\cal F}(x,t,Q)+1.
\end{equation}
For simplicity, in
what follows we may write ${\cal G}(Q)$ in place of ${\cal
G}(x,t,Q)$. For a given function $Q_0(x,t)$ on $\mathbb{R}^+ \times
\mathbb{R}^+$, we study the sequence of iterates of ${\cal G}$ given
by \begin{equation} \label{defQm} Q_{m+1}(x,t)={\cal G}(x,t,Q_m)
\quad \mbox{for $m \ge 0$}. \end{equation} The main result of this
section, given in Theorems \ref{thm:Pcmt} and \ref{fde-sol-unique}, proves the existence and uniqueness of a
solution~$\tQ$ to the equation \begin{equation}
\label{tQ-fixed-point} \tQ(x,t)={\cal G}(x,t,\tQ),
\end{equation}
and the
exponential rate of convergence of the sequence $Q_m$ to $\tQ$ when $Q_0$ is continuous and $a$
and $\otimes$ satisfy certain continuity-type properties. These conditions make ${\cal G}$
into a type of contraction, allowing us to apply the methods of the
contraction mapping theorem (see, for instance, \citet{Luenberger69}) to
prove that sequences of functions generated by the iterates of ${\cal G}$
on some initial function converge to a unique fixed point of ${\cal G}$.

For a function $Q$ bounded over all compact domains  $\D \subset
\mathbb{R^+} \times \mathbb{R^+}$, let
$$
||Q||_\D=\sup_{(y,s) \in \D} |Q(y,s)|.
$$
For any two functions $Q$ and $R$  defined on $\D \subset
\mathbb{R}^+ \times \mathbb{R}^+$ and $\delta \in [0,\infty)$, let
$${\bf d}_{\D,\,\delta}^1(Q,R)=\sup_{\{(x,t),(x',t) \in \D\,
:\, |x-x'|\le \delta\}}|Q(x,t)-R(x',t)|,
$$
and
$${\bf d}_{\D,\,\delta}(Q,R)=\sup_{\{(x,t),(x',t') \in \D\,
:\, |x-x'|\le \delta, |t-t'| \le \delta \}}|Q(x,t)-R(x',t')|.
$$
The functions ${\bf d}_{\D,\,\delta}^1$ and ${\bf d}_{\D,\,\delta}$ satisfy all the properties of a metric except that they are not necessarily zero when $Q=R$.

\begin{definition}
Given a function $a(y)$, the operation $\otimes$ is of \emph{$\kappa$-contraction type} for some $\kappa\ge 0$ if, for all compact
domains $\D = [0,X] \times [0,T]$,
there exist $\beta \ge 0$ and
$\eta(\delta)$, defined for $\delta \ge 0$ and satisfying
\beas
\eta(0)=0 \quad \mbox{and} \quad \lim_{\delta \downarrow 0} \eta(\delta)=0,
\enas
such that for
any functions $Q$ and $R$ on $\mathbb{R}^+ \times \mathbb{R}^+$,
\begin{equation} \label{eq:cont3}
{\bf d}_{\D,\,\delta}^1((a \otimes Q),(a \otimes R))\le \kappa {\bf
d}_{\D,\,\delta}^1(Q,R)+\eta(\delta)(||Q||_D \vee ||R||_D+\beta),
\end{equation}
for all $\delta \in [0,\infty)$, where $b \vee c$ is the maximum of $b$ and $c$.
The operation $\otimes$ is
of \emph{contraction type} if it is of $\kappa$-contraction type for some
$\kappa\ge 0$.
\end{definition}

\begin{lemma}
\label{QRsup} Let ${\cal G}$ be given by (\ref{def:calG}) or (\ref{def:calFG}) for
$\otimes$ of $\kappa$-contraction type. Let $X \ge 0$, $T' \ge T \ge
0$ be arbitrary, and set $\D\,'=[0,X] \times [0,T']$ and $\D=[0,X]
\times [0,T]$. If $Q$ and $R$ are any functions defined on
$\D{\,}'$, then
\begin{equation} \label{eq:cont1} \vert \vert {\cal G}(Q)-{\cal G}(R) \vert \vert_{\D'}
\le \kappa T\vert \vert Q-R \vert \vert_{\D} + \kappa(T'-T)\vert
\vert Q - R \vert \vert_{\D\,'}, \end{equation} and when $a(y)$ is bounded, for all $\delta
\ge 0$, the inequality
\begin{equation}
\label{eq:cont2} {\bf d}_{\D,\,\delta}({\cal G}(Q),{\cal G}(R))
 \le
\kappa T{\bf d}_{\D,\,\delta}^1(Q,R) +\omega(\delta)(||Q||_D \vee ||R||_D+\gamma) \end{equation}
holds for some $\gamma \ge 0$ and $\omega(\delta)$ satisfying $\omega(0)=0$ and $\lim_{\delta \downarrow 0}\omega(\delta)=0$.
\end{lemma}

\proof Clearly it suffices to prove the claim for ${\cal G}$ of the form (\ref{def:calG}).
To obtain \eqref{eq:cont1} let $(x,t) \in \D\,'$. Applying
\eqref{eq:cont3} with $\delta=0$, and noting ${\bf d}_{\D,0}^1( Q,
R)=\vert \vert Q-R \vert \vert_{\D}$, we have \beas
\lefteqn{\vert {\cal G}(x,t,Q)-{\cal G}(x,t,R) \vert} \\
&\le&  \int_0^t \vert (a \otimes Q)(x,s)-(a \otimes R)(x,s) \vert ds
\le \int_0^{T'} \vert (a \otimes Q)(x,s)-(a \otimes R)(x,s) \vert ds\\
&=&  \int_0^T \vert (a \otimes Q)(x,s)-(a \otimes R)(x,s) \vert ds +
\int_T^{T'} \vert (a \otimes Q)(x,s)-(a \otimes R)(x,s) \vert ds \\
&\le&  \int_0^T \kappa {\bf d}_{\D,0}^1( Q, R) ds + \int_T^{T'}
\kappa {\bf d}_{\D',0}^1( Q, R)  ds\\
&\le&  \kappa T\vert \vert Q-R \vert \vert_{\D} + \kappa (T'-T)\vert
\vert Q  - R \vert \vert_{\D\,'}. \enas Taking supremum on the left
hand side over $(x,t)\in \D\,'$, we obtain \eqref{eq:cont1}.

To obtain \eqref{eq:cont2}, let $\delta \ge 0$ and suppose $(x,t),
(x',t') \in \D$ satisfy $|x-x'| \le \delta$ and $|t-t'| \le \delta$;
assume first that $t' \ge t$. Then, applying \eqref{eq:cont3}, for each of the four cases of
\eqref{twostars} we
have \beas
\lefteqn{\vert {\cal G}(x',t',Q)-{\cal G}(x,t,R) \vert} \\
&\le&  \int_0^t \vert (a \otimes Q)(x',s)-(a \otimes R)(x,s) \vert ds + \int_t^{t'}\vert (a \otimes Q)(x',s)\vert ds\\
&\le&  T \left( \kappa {\bf d}_{\D,\,\delta}^1(Q,R) +  \eta(\delta)(||Q||_D \vee ||R||_D)+\beta \right) + \delta ||a|| (e^T+  \vert \vert Q
\vert \vert_{\D}), \enas where $||a||$ is the supremum of $a$ on $[0,\infty)$. If we assume that $t' \le t$ we obtain a
similar bound, except that the last term has $\vert \vert R \vert
\vert_{\D}$ in place of $\vert \vert Q \vert \vert_{\D}$. Taking the
above two cases into account and taking supremum on the left hand
side, we obtain \eqref{eq:cont2} with $\gamma=\beta+ e^T$ and $\omega(\delta)=T\eta(\delta)+\delta ||a||$. \qed

\begin{definition}
\label{Gkaptype}
If ${\cal G}$ satisfies the conclusions of Lemma \ref{QRsup} for some $\kappa \ge 0$, $\gamma \ge 0$, and $\omega(\delta)\ge 0$, then ${\cal G}$ is of \emph{$\kappa$-contraction type}. The operator ${\cal G}$ is of \emph{contraction type} if it is of $\kappa$-contraction type for some $\kappa \ge 0$.
\end{definition}

Tautologically, under Definition \ref{Gkaptype},
Lemma \ref{QRsup} says that if ${\cal G}$ is given by (\ref{def:calG}) with $\otimes$ of contraction type, then ${\cal G}$ is of contraction type.

\begin{proposition}
\label{vectorspace}
If ${\cal G}_1$ and ${\cal G}_2$ are of
$\kappa_1$ and $\kappa_2$-contraction type, respectively, then ${\cal G}_1+{\cal G}_2$ is of $(\kappa_1+\kappa_2)$-contraction  type, and $\alpha {\cal G}_1$ is of $|\alpha| \kappa_1$-contraction type for any $\alpha \in \mathbb{R}$. Hence the collection of contraction type operators form a vector space over $\mathbb{R}$. The constant operator ${\cal G}(\cdot)=1$ is a contraction operator of type $0$.
\end{proposition}

\proof
The triangle inequality shows the first claim, the remaining claims are self evident.
\qed

We note that Proposition \ref{vectorspace} shows that if ${\cal G}$ is of contraction type then so is ${\cal G}+1$.
The next corollary shows that any ${\cal G}$ of contraction type preserves continuity.
\begin{corollary}
\label{recursion-continuous} If  $a(y)$ is bounded, $Q$ is continuous on $\mathbb{R}^+ \times \mathbb{R}^+$, and ${\cal G}$ is of contraction type, then ${\cal G}(Q)$ is continuous.
\end{corollary}
\proof Given $\varepsilon>0$, by the uniform continuity of $Q$ in
any compact domain $\D$ of the form $[0,X] \times [0,T]$, and, since $\omega(\delta) \downarrow 0$ as $\delta \downarrow 0$, there exists $\delta>0$
such that $\kappa T{\bf d}_{\D,\,\delta}^1(Q,Q)<\varepsilon/2$ and
$\omega(\delta)(||Q||_D + \gamma)  <\varepsilon/2$. Hence
${\bf d}_{\D,\,\delta}({\cal G}(Q),{\cal G}(Q))<\varepsilon$ by
(\ref{eq:cont2}). \qed

\begin{theorem}
\label{thm:Pcmt} Let $Q_0(x,t)$ be any continuous function on $\mathbb{R}^+
\times \mathbb{R}^+$, and let $Q_m, m=1,2,\ldots$, be given by
(\ref{defQm}) for some ${\cal G}$ of $\kappa$-contraction type. Then $Q_m$, $m=0,1,\ldots$,
are continuous and converge uniformly on compact subsets of
$\mathbb{R}^+ \times \mathbb{R}^+$ to a continuous function $\tQ$ that
satisfies (\ref{tQ-fixed-point}).

Additionally, for every compact domain $\D \subset \mathbb{R}^+
\times \mathbb{R}^+$ there exists a constant $C>0$ and $\theta \in
(0,1)$ such that
$$
||\tQ-Q_m||_\D \le C\theta^m \quad \mbox{for all $m=0,1,\ldots$.}
$$
\end{theorem}

\proof Let $\kappa>0$ as otherwise $Q_m=0$ for all $m \ge 1$ and the
result is trivial. Let $X \ge 0$ be arbitrary, and let ${\widetilde
T}$ be the supremum over all $T \ge 0$ such that there exist $C>0$
and $\theta \in (0,1)$ such that
\begin{equation} \label{Cthetam} c_{m,T} \le C\theta^{m} \quad \mbox{for
$m=1,2,\ldots,$ where} \quad c_{m,T}=\vert \vert Q_m-Q_{m-1} \vert
\vert_{{[0,X] \times [0,T]}}\,.
\end{equation} We first prove that ${\widetilde T}=\infty$. Let $\theta \in (0,1)$. Using
\eqref{eq:cont1} of Lemma \ref{QRsup} with $T=0$ and $T'=\theta/\kappa$
we obtain $$ c_{m,T'} \le \left( \frac{c_{1,T'}}{\theta}\right)
\theta^m \quad \mbox{for $m=1,2,\ldots$} $$ by induction, first noting
that it holds trivially for $m=1$.  Hence ${\widetilde T}>0$.

If (\ref{Cthetam}) holds for some ${\cal C }>0$, $\vartheta \in
(0,1)$ and $T>0$ then, for $T'=T+\vartheta/\kappa$, using
\eqref{eq:cont1} again \begin{equation} \label{inducta} c_{m,T'} \le
\kappa Tc_{m-1,T}+ \kappa (T'-T)c_{m-1,T'} \le \kappa T{\cal C
}\vartheta^{m-1} + \vartheta c_{m-1,T'} \quad \text{for} \quad m \ge
2.\end{equation} The inequality $$ c_{m,T'} \le C' m \vartheta^{m-1}
\quad \mbox{where} \quad C'=\max\{\kappa T{\cal C},c_{1,T'}\} $$
holds trivially for $m=1$, and now induction, using (\ref{inducta}),
shows it holds for all $m \ge 1$. Elementary calculus now shows that
$c_{m,T'} \le C' m \vartheta^{m-1}$ for $m\ge 1$ implies that for
any $\theta \in (\vartheta,1)$ there exists $C$ such that
$c_{m,T'} \le C \theta^m$ for all $m \ge 1$. Hence \eqref{Cthetam}
holds  when replacing $T$ by $T'$, and therefore assuming
${\widetilde T}<\infty$ leads to a contradiction. Now note that
\eqref{Cthetam} holds when replacing $[0,X] \times [0,T]$ by any
compact $\D$, and $c_{m,T}$ by $c_{m,\D}=||Q_m-Q_{m-1}||_\D$, as $D \subset [0,X] \times [0,T]$ for some $X$ and $T$
sufficiently large.

By Corollary \ref{recursion-continuous} and induction, the functions
$Q_m$, $m=0,1,\ldots$, are continuous and by (\ref{Cthetam}) they form a Cauchy sequence in the supremum norm on any compact domain $\D$. Therefore $Q_m$, $m=0,1,\ldots$, converges uniformly on any compact domain to
a continuous limit function $\tQ$. The function $\tQ$ solves
(\ref{tQ-fixed-point}) as applying the triangle inequality and
then \eqref{eq:cont1} of Lemma \ref{QRsup} with $T'=T$ yields \beas
||{\cal G}(\tQ)-\tQ||_\D &\le& ||{\cal G}(\tQ)-Q_{m+1}||_\D+||Q_{m+1}-\tQ||_\D\\
&=& ||{\cal G}(\tQ)-{\cal G}(Q_m)||_\D+||Q_{m+1}-\tQ||_\D\\
&\le& \kappa T||\tQ-Q_m||_\D+||Q_{m+1}-\tQ||_\D, \enas whose right
hand side converges to zero as $m \rightarrow \infty$. Hence ${\cal G}(\tQ)=\tQ$.

Lastly, by the triangle inequality, for all $n \ge m$, by (\ref{Cthetam}) $$
||Q_n-Q_m||_\D \le C\sum_{j=m+1}^n \theta^{j} \le \left(
\frac{C}{1-\theta} \right) \theta^{m+1},
$$
so letting $n \rightarrow \infty$ yields the final claim of the theorem.
\qed

\begin{theorem}
\label{fde-sol-unique} If ${\cal G}$ is of contraction type, equation (\ref{tQ-fixed-point}) has a unique
solution that is bounded over all compact domains of $\mathbb{R^+} \times
\mathbb{R^+}$.
\end{theorem}

\proof Let $Q$ and $R$ be two solutions to (\ref{tQ-fixed-point})
that are bounded over all compact subsets of $\mathbb{R^+} \times
\mathbb{R^+}$. Let $\D=[0,X] \times [0,T]$ and $\D'=[0,X]
\times [0,T']$, for some $X \ge 0 $ and $T' \ge T \ge 0$. Taking $T=0$ and $\kappa T'<1$ in
\eqref{eq:cont1} yields $||Q-R||_{\D'}=0$.

Replacing $T$ and  $T'$  by $T'$ and $T''$, respectively, in
\eqref{eq:cont1}, where $0 \le \kappa(T''-T') < 1$, it is easy to see that for the
larger domain $\D''=[0,X] \times [0,T'']$ we now obtain
\begin{equation}\label{eq:uni} ||Q-R||_{\D''} = ||{\cal G}(Q)-{\cal G}(R)||_{\D''} \le \kappa (T''-T') ||Q-R||_{\D''}. \end{equation}
This implies that $||Q-R||_{\D''}=0$, and it
follows that $||Q-R||_{\D}=0$ holds for any domain $[0,X]
\times [0,T]$, and hence for any compact domain $\D$. \qed

In the  cases of interest to us here, the operator ${\cal G}$ has
some additional properties, in particular, monotonicity
and a type of boundedness.

\begin{definition}\label{def:bdmon} An operator ${\cal G}$ is \emph{$B$-bounded} for a continuous function $B:
\mathbb{R}^+ \times \mathbb{R}^+ \rightarrow \mathbb{R}^+$ if $$ Q
\le B \quad \mbox{implies} \quad {\cal G}(Q) \le B. $$
An operator ${\cal G}$ is \emph{monotone} if $$ Q \le R \quad \mbox{implies}
\quad {\cal G}(Q) \le {\cal G}(R). $$
\end{definition}
Clearly the operator ${\cal G}$ of \eqref{def:calG} is monotone for
all four choices in \eqref{twostars}; hence so is ${\cal G}$ of \eqref{def:calFG} with
$\otimes$ given by case $(Bo)$. The following proposition is
immediate from the preceding definitions.

\begin{proposition}
\label{sandwich}
Let ${\cal G}$ be a monotone, $B$-bounded operator of contraction type, and let $Q_m$, $m=0,1,\ldots$,
be the functions given by (\ref{defQm})
when $Q_0(x,t)\equiv0$, and $Q_m'$ the functions
when $Q_0'(x,t)=B(x,t)$. Then the unique solution $\tQ$ to (\ref{tQ-fixed-point}) satisfies
$$
0 \le Q_m(x,t) \le \tQ(x,t) \le Q'_m(x,t) \le B(x,t) \quad \mbox{for all $m=0,1,\ldots$.}
$$
\end{proposition}

If the iteration (\ref{defQm}) is used for  numerical
approximation of the solution $\tQ$ to (\ref{tQ-fixed-point}), then Proposition
\ref{sandwich} allows for the verification of the precision of the calculation, as it
provides a lower and upper bound for $\tQ$; the distance between these
bounds therefore gives the approximation error.

We now show that when $a$ is a bounded function, ${\cal G}$ of \eqref{def:calG} is of contraction type in
all cases of \eqref{twostars}, and that the same is true for ${\cal G}$ of \eqref{def:calFG} with $\otimes$
given in case $(Bo)$.

\ignore{For the first case, that of the Bomber, write \begin{equation}
\label{def-calF} {\cal G}(x,t,Q)=\int_0^t (a \star Q)(x,s) ds \quad
\mbox{and} \quad {\cal F}(x,t,Q)=\int_0^t (a \star Q)(x,s) ds+1, \end{equation}
where $\star$ is given by the max convolution $$ (a \star Q)(x,t) =
\sup_{0 \le y \le x}a(y)Q(x-y,t). $$ }

\begin{lemma}
\label{lem:F-is-ctype} Let $a: \mathbb{R}^+ \rightarrow \mathbb{R}^+$
be uniformly continuous
and bounded by $\kappa$. The operator ${\cal G}$ given by (\ref{def:calG}) is of
$\kappa$-contraction type
in cases $(Bo)$ and $(F_0)$ of
(\ref{twostars}), of 1-contraction type in case $(F_1)$, and of $[(1-u)\kappa +u]$-contraction type in case $(F_u)$.
The operator ${\cal G}$ given by (\ref{def:calFG}) with $\otimes$ as in case $(Bo)$ is of $\kappa$ contraction type.
\end{lemma}

\proof 
For the first claim, with $\otimes$ given by case $(Bo)$ of \eqref{twostars} it suffices, by Lemma \ref{QRsup}, to prove
that $\otimes$ is of $\kappa$-contraction type. Let a compact $D =[0,X] \times [0,T]$ and $\delta \ge 0$ be given and $(x,t),(x',t) \in D$ with $|x-x'| \le \delta$.  For a given $\epsilon>0$ let
$y^* \in [0,x']$ achieve $\sup_{0 \le y \le x'}a(y)Q(x'-y,t)$ to within $\epsilon$, and set  $y_1^*=\min\{y^*,x\}$. Hence
\beas
\lefteqn{(a \otimes Q)(x',t)-(a \otimes R)(x,t)}\\
&\le& a(y^*)Q(x'-y^*,t)-a(y_1^*)R(x-y_1^*,t)+\epsilon \\
&=& a(y^*)[Q(x'-y^*,t)-R(x-y_1^*,t)]+(a(y^*)-a(y_1^*))R(x-y_1^*,t) + \epsilon.
\enas
Since $0 \le a(y) \le \kappa$ for all $y \ge 0$ and $|(x'-y^*)-(x-y_1^*)| \le \delta$, with
\beas
\eta(\delta)=\sup_{|y'-y|\le \delta}|a(y')-a(y)|,
\enas
we obtain
\bea \label{aotimesQeta}
(a \otimes Q)(x',t)-(a \otimes R)(x,t) \le \kappa {\bf d}_{\D,\delta}^1(Q,R) + \eta(\delta)||R||_D.
\ena

Similarly, there exists another
$y^* \in [0,x]$ and corresponding $y_1^* \in [0,x']$ such that
$$ (a \otimes Q)(x',t)-(a \otimes R)(x,t) \ge
a(y_1^*)Q(x'-y_1^*,t)-a(y^*)R(x-y^*,t)]-\epsilon,
$$
which yields a bound of the form (\ref{aotimesQeta}) with the roles of $Q$ and $R$ interchanged.
Since $\epsilon$ is
arbitrary
$$
|(a \otimes Q)(x',t)-(a \otimes R)(x,t)| \le \kappa {\bf d}_{\D,\delta}^1(Q,R) + \eta(\delta)(||Q||_D \vee ||R||_D).
$$
Taking supremum on the left hand side and noting that $\eta(0)=0$ and $\lim_{\delta \downarrow 0} \eta(\delta)=0$ by the uniform continuity of $a(y)$ yields the desired conclusion. The first and final claims of
the lemma are now immediate.

In case $(F_u)$ we have in a similar way
\beas
\lefteqn{(a \otimes Q)(x',t)-(a \otimes
R)(x,t)}\\
&\le& a(y^*)e^t+[a(y^*)+u(1-a(y^*))]Q(x-y^*,t)- a(y_1^*)e^t-[a(y_1^*)+u(1-a(y_1^*))]R(x-y_1^*,t) +\epsilon.
\enas
Adding and subtracting terms as before, we obtain
\beas
&& (a(y^*)-a(y_1^*))e^t+[(1-u)a(y^*)+u][Q(x-y^*,t)-R(x-y_1^*,t)]\\
&+&[(a(y^*)-a(y_1^*))(1-u)]R(x-y_1^*,t) +\epsilon\\
&\le& [(1-u)\kappa +u] {\bf d}_{\D,\delta}^1(Q,R)+ \eta(\delta)e^T(||R||_D+1) + \epsilon,
\enas
since $|(1-u)a(y^*)+u| \le (1-u)\kappa +u$. The result
follows as for case $(Bo)$.  Now cases $(F_0)$ and $(F_1)$ follow by specializing $(F_u)$
to $u=0$ and $u=1$, respectively. \qed

\begin{lemma}
\label{R-is-bounded} Let ${\cal G}$ be given by (\ref{def:calFG}) for
case $(Bo)$ of (\ref{twostars}).
If $a: \mathbb{R}^+ \rightarrow \mathbb{R}^+$ is bounded by $\kappa$ then ${\cal G}$ is $e^{\kappa t}$-bounded.
\end{lemma}

\proof  If $Q(x,s) \le e^{\kappa s}$, using that $0 \le a(y) \le \kappa$ for
all $y \ge 0$, we have for all $x \ge 0$,
 $$ {\cal G}(x,t,Q) \le
\kappa \int_0^t \sup_{0 \le y \le x} Q(y,s)ds + 1 \le  \kappa \int_0^t e^{\kappa s} ds + 1 =
e^{\kappa t}, $$ proving the boundedness claim.
\qed

\ignore{ 
\begin{lemma}
\label{DQRsup} The operator ${\cal E}$ given by (\ref{defcalE}) is of
$1$-contraction type.
\end{lemma}

\proof By Lemma \ref{QRsup}, if suffices to show that the operation $*$ given by (\ref{twostars}) is of $1$-contraction type.
The proof is the essentially the same as that of Lemma \ref{QRsup}. Let $\delta \ge 0$ and $(x,t),(x',t)\in \mathbb{R^+} \times \mathbb{R^+}$ with $|x-x'| \le \delta$, and again let $Q(x,t)=Q(x \vee 0,t)$.  For every $t$ and every $\epsilon>0$ there exists
$y^* \in [0,x']$, nearly achieving the supremum $(a*Q)(x,t)$,
such that $$ (a * Q)(x',t)-(a * R)(x,t) \le
[Q(x'-y^*,t)-R(x-y^*,t)]+\epsilon. $$ Similarly, there exists
$y^* \in [0,x]$ nearly achieving the supremum $(a*R)(x,t)$,
and for which the inequality is reversed, now yielding the conclusion just
as in the proof of Lemma \ref{QRsup}.
\qed

} 

\begin{lemma}
\label{E-is-bounded}
Consider case $(F_u)$ of \eqref{twostars} with $a:
\mathbb{R}^+ \rightarrow \mathbb{R}^+$ bounded by $\kappa$, and let
$\alpha=(1-u)\kappa+u$.
Then
${\cal G}$ of (\ref{def:calG}) is $B$-bounded where $B(x,t)=B_\alpha(t)$ is given by
\beas
B_\alpha(t)=\left\{
\begin{array}{cc}
\kappa t e^t & \alpha=1\\
\frac{\kappa }{1-\alpha}\left(e^t-e^{\alpha t}\right) & \alpha \not =1.
\end{array}
\right.
\enas
In particular, in case $(F_0)$ the operator ${\cal G}$ is $B_\kappa(t)$-bounded,
and in case $(F_1)$ is $\kappa te^t$-bounded.
\end{lemma}

\proof It is easily checked in both the $\alpha=1$ and $\alpha \not =1$ cases that if $Q(x,t) \le B_\alpha(t)$, then
$$
{\cal G}(x,t,Q) \le \int_0^t \left( \kappa e^s + \alpha B_\alpha(s) \right) ds = B_\alpha(t).
$$
\qed

The following corollary summarizes the results above.
\begin{corollary}
\label{sec:final-corollary} Let $a:\mathbb{R}^+ \rightarrow
\mathbb{R}^+$ be uniformly continuous and bounded by $\kappa$. Then, with ${\cal G}$ given
by (\ref{def:calG}) and $\otimes$ any of the four cases of (\ref{twostars}),
or ${\cal G}$ given by (\ref{def:calFG}) with $\otimes$ given in case $(Bo)$ of (\ref{twostars}),
there exists a continuous
solution to the equation ${\cal G}(Q)=Q$, and there are no other solutions bounded over compact
domains. Moreover, with $Q_0(x,t)$ any continuous function from $\mathbb{R}^+ \times \mathbb{R}^+$
to $\mathbb{R}^+$, the
sequence of iterates of ${\cal G}$ given by $Q_{m+1}={\cal G}(Q_m)$
for $m \ge 0$ converges uniformly to the solution of ${\cal G}(Q)=Q$
at an exponential rate. For every $m$, the solution to ${\cal G}(Q)=Q$
is bounded below
by the $m^{th}$ term of the sequence initialized with $Q_0(x,t)=0$,
and if ${\cal G}$ is $B$-bounded, then it is also
bounded above by the $m^{th}$ term of the sequence initialized
with $Q_0(x,t)=B(x,t)$.
\end{corollary}

We conclude this section with a proof of Theorem \ref{thm:genFighter}.

\noindent \textit{Proof of Theorem \ref{thm:genFighter}.} Let ${\cal G}$ be given by (\ref{def:calG}) and $\otimes$ by $(F_u)$ of (\ref{twostars}) for $u \in [0,1]$. By Lemma \ref{lem:F-is-ctype}, when $a$ is a probability, the operator ${\cal G}$ is $te^t$-bounded. By Corollary \ref{sec:final-corollary} the unique solution $\tN$ to ${\cal G}(Q)=Q$ exists, is continuous and satisfies $0 \le \tN(x,t) \le te^t$.
Now letting $\tN(x,t)=e^tN(x,t)$ as in (\ref{tP=etP}), we have $N(x,t) \in [0,t]$, as claimed. \qed

\section{The Bomber}
\label{sec:bomber} In this section we consider the dynamic
programming equation (\ref{tP-integral-equation}), and
specialize the results of the previous section to $\otimes$ given by
case $(Bo)$ of (\ref{twostars}). Let ${\cal F}(Q)$ and
${\cal G}(Q)$ be given by
(\ref{def:calFG}). In
particular, equation (\ref{tP-integral-equation}) may be written
more compactly as ${\cal G}(\tP)=\tP$. Throughout this section
$a(y)$ can be interpreted as the probability of the Bomber surviving
an enemy encounter after expending an amount~$y$ of its ammunition, although our results hold in greater generality, in particular,
 for any bounded log-concave $a:\mathbb{R}^+\To
\mathbb{R}^+$. We require strict log-concavity for
Theorem~\ref{C-dcb}. For example, the function $a(y)$ in (\ref{a-function})
is strictly log-concave in $y>0$ for any $u<1$, being the composition of a strictly concave function
and a strictly increasing and concave function, the log.


We show that ${\cal F}$ and ${\cal G}$ preserve $\mbox{TP}_2$ and
other properties, which are used to prove Theorem~\ref{thm:A}, that
[A] holds for the  Bomber problem. The Fighter Problems that we
treat in subsequent sections have much in common with the present
one. We begin by collecting some facts about $\mbox{TP}_2$ and log
concave functions in the following lemma; the claims which cannot be
easily verified using definition~(\ref{TP2-def}) can be found in
\citet{Karlin68}. For part~1 of Lemma~\ref{lem:tp2} see also~\citet{Schoenberg51}.
Note also that part~2 is a special case of part~3, and that the second conclusion of part~5 follows from part~4 after noting that a function $f(x)$ of the single variable~$x$, or $g(y)$ of $y$, is $\mbox{TP}_2(x,y)$.
\begin{lemma}
\label{lem:tp2}

\begin{enumerate}

\item A nonnegative function $a(y)$ is
(strictly) log-concave if and only if $a(x-y)$ is (strictly) $\mbox{TP}_2(x,y)$.

\item The convolution of log-concave functions is log-concave.

\item (The composition formula) If $J(x,y)$ is $\mbox{TP}_2(x,y)$,
$L(y,z)$ is $\mbox{TP}_2(y,z)$, $\sigma$ is a
nonnegative sigma-finite measure on $\mathbb{R}$, and $M(x,z)=\int
J(x,y)L(y,z)d\sigma(y)$ exists, then $M(x,z)$ is $\mbox{TP}_2(x,z)$.

\item If $J(x,y)$ and $L(x,y)$ are
$\mbox{TP}_2(x,y)$, then so is their product.

\item If $J(x,y)$ is
$\mbox{TP}_2(x,y)$ then $J(x,y)^p$ for any positive $p$, and
$J(x,y)f(x)g(y)$ for any nonnegative functions $f$ and $g$, are $\mbox{TP}_2$.
\end{enumerate}
\end{lemma}

\ignore{For the given function $a(y)$ on $\mathbb{R}^+$, we define $a(y)=0$ for $y<0$; if $a$ is log-concave
on $\mathbb{R}^+$ the extension is log-concave on all of
$\mathbb{R}$. Similarly }
Given a $\mbox{TP}_2$ function $Q$ defined on $\mathbb{R}^+ \times \mathbb{R}^+$ we extend it to
$\mathbb{R}  \times \mathbb{R}$ by setting it equal to zero where it was not defined before; the resulting
function is $\mbox{TP}_2$ on $\mathbb{R} \times \mathbb{R}$. The next lemma follows directly from \eqref{def:calG} and \eqref{def:calFG}.

\begin{lemma}
\label{R-is-increasing} If $a: \mathbb{R}^+ \rightarrow \mathbb{R}^+$ and $Q(x,t)$ is
nonnegative and increasing in $x$ then ${\cal F}(Q)$ and ${\cal
G}(Q)={\cal F}(Q)+1$ are nonnegative and increasing in $x$ and
$t$.
\end{lemma}
\ignore{ BEGIN IGNORE
\proof The claims on nonnegativity and monotonicity follows directly from \eqref{def:calG} and \eqref{def:calFG}.
If $Q$ is increasing in $x$, then
for $x' \ge x$ we have
$$
\fmx_{0 \le y \le x'}a(y)Q(x'-y,s) \ge \fmx_{0 \le y \le x}a(y)Q(x'-y,s) \ge \fmx_{0 \le y \le x}a(y)Q(x-y,s),
$$
and therefore, with $t' \ge t$, because $a$ and $Q$ are nonnegative,
\beas
{\cal F}(x',t',Q) &=& \int_0^{t'} \fmx_{0 \le y \le x'} a(y)Q(x'-y,s)ds\\
                &\ge& \int_0^t   \fmx_{0 \le y \le x} a(y)Q(x-y,s)ds={\cal F}(x,t,Q).
\enas Hence ${\cal G}(x,t,Q) = {\cal F}(x,t,Q)+1$ is also
increasing in $x$ and $t$ . \qed  END IGNORE }

The following key lemma gives conditions under which the $\mbox{TP}_2$ property is preserved by the sup-convolution. We recall that $\otimes$ in this section is given by $(Bo)$ of (\ref{twostars}).
\begin{lemma}
\label{maxconvpreserves} If $Q$ is a continuous, nonnegative
$\mbox{TP}_2$ function on $\mathbb{R}^+ \times \mathbb{R}^+$ and $a:
\mathbb{R}^+ \rightarrow \mathbb{R}^+$ is log-concave, then $a \otimes
Q$ is $\mbox{TP}_2$.
\end{lemma}

\proof Let $(x,t) \in \mathbb{R}^+ \times \mathbb{R}^+$ and $p>0$ be arbitrary. Making the change of variable $v=x-y$ and writing the limits of integration using
the indicator $I(x,v)={\bf 1}(0 \le v \le x)$ we obtain
$$
\int_0^x \left[ a(y)Q(x-y,s) \right]^p dy = \int_{-\infty}^\infty  a(x-v)^p I(x,v)
 Q(v,s)^p dv.
$$
One can verify directly from (\ref{TP2-def}) that $I(x,v)$ is $\mbox{TP}_2(x,v)$. Since $a(x)$ is log-concave on $\mathbb{R}^+$, defining $a(x)=0$ for $x<0$ we have $a(x)$ log-concave on all of $\mathbb{R}$, and so
the product $a(x-v)I(x,v)$, and therefore also $a(x-v)^pI(x,v)$, are $\mbox{TP}_2(x,v)$ on $\mathbb{R} \times \mathbb{R}$.
Hence the integrand is a product of a
$\mbox{TP}_2(x,v)$ function with a $\mbox{TP}_2(v,s)$ function, and
so the integral is $\mbox{TP}_2(x,s)$ by part 3 (the composition formula) of
Lemma \ref{lem:tp2}. Raising to the power $1/p$, we therefore have
\begin{equation} \label{Lpnorms} \left( \int_0^x \left[ a(y)Q(x-y,s)
\right]^p dy \right)^{1/p} \quad \mbox{is $\mbox{TP}_2(x,s)$ for all
$p>0$.} \end{equation}

Since $a(y)$ is log-concave on $\mathbb{R}^+$ it is continuous in $(0,x]$ and bounded on $[0,x]$, with a possible
discontinuity at $0$, in which case there is a downward jump. Hence for all fixed
$x$ and $s$ the product $a(y)Q(x-y,s)$ has the same properties. Therefore the
limit of the $L^p$ norms (\ref{Lpnorms}) of the function
$a(y)Q(x-y,s)$ of $y$ on $[0,x]$ converges to its essential
supremum which coincides with its maximum value over $y$. Using the
nonnegativity of $a(\cdot)$ and $Q(\cdot,\cdot)$, we therefore have,
by the preservation of the $\mbox{TP}_2$ property under limits, that
$$
\lim_{p \rightarrow \infty}\left( \int_0^x \left[
a(y)Q(x-y,s) \right]^p dy \right)^{1/p} = (a \otimes Q)(x,s) \quad \mbox{is $\mbox{TP}_2(x,s)$.}
$$
\qed

\begin{lemma}
\label{integral-TP2} If $a: \mathbb{R}^+ \rightarrow \mathbb{R}^+$
is log-concave and $Q$ is $\mbox{TP}_2$ then ${\cal F}(Q)$ is
$\mbox{TP}_2$.
\end{lemma}
\proof By Lemma \ref{maxconvpreserves} $(a \otimes Q)(x,s)$ is
$\mbox{TP}_2(x,s)$. As $I(s,t)={\bf 1}(0 \le s \le t)$ is
$\mbox{TP}_2(s,t)$, the integral of the product
$$
{\cal F}(x,t,Q)=\int_0^t  (a \otimes Q)(x,s)ds = \int_{-\infty}^\infty  (a \otimes Q)(x,s)I(s,t)ds,
$$
is $\mbox{TP}_2(x,t)$.
\qed

\begin{lemma}
\label{plus1-TP2} If $Q(x,t)$ is nonnegative, increasing
in $x$ and $t$, and $\mbox{TP}_2$, then so is $Q(x,t)+1$.
\end{lemma}

\proof It is clear that $Q(x,t)+1$ is nonnegative and increasing in $x$ and $t$
whenever $Q(x,t)$ is. To demonstrate the $\mbox{TP}_2$ property it is required to show that
$$
(Q(x',t')+1)(Q(x,t)+1) \ge  (Q(x',t)+1)(Q(x,t')+1),
$$
where $x'>x$ and $t'>t$. The monotonicity and nonnegativity of $Q$ imply that the result is
trivial if $Q(x',t')=0$. Otherwise using $Q(x',t')Q(x,t)\ge Q(x',t)Q(x,t')
$ we obtain \beas
\lefteqn{(Q(x',t')+1)(Q(x,t)+1) - (Q(x',t)+1)(Q(x,t')+1)}\\
&\ge& (Q(x',t')+1)(\frac{Q(x',t)Q(x,t')}{Q(x',t')}+1) - (Q(x',t)+1)(Q(x,t')+1)\\
&=&(Q(x',t')-Q(x',t))(Q(x',t')-Q(x,t'))/Q(x',t') \ge 0
\enas
where the last inequality
follows since
$$
Q(x',t') \ge \max\{ Q(x',t),Q(x,t')\},
$$
as $Q$ is increasing in both $x$ and $t$.\qed

\ignore{By cancelling common factors, and using
\begin{equation} \label{def-A}
Q(x',t')Q(x,t)\ge  Q(x',t)Q(x,t')
\end{equation}
we see that it suffices to show
\begin{equation} \label{A-sums}
Q(x',t')+Q(x,t)\ge  Q(x',t)+Q(x,t').
\end{equation}
Letting $B(x,t)=\log Q(x,t)$, inequality (\ref{def-A}) yields
$$
B(x',t')+B(x,t) \ge  B(x',t)+B(x,t').
$$
By the monotonicity property of $Q(x,t)$, we have that $B(x',t')$
dominates both $B(x',t)$ and $B(x,t')$. Clearly there exists $b
\le B(x,t)$ such that
$$
B(x',t') + b = B(x',t)+B(x,t').
$$
Hence, the vector $(B(x',t'),b)$ dominates $(B(x',t),B(x,t'))$ in
the majorization order, and therefore, for any convex function $f$, in
particular for $f(u)=e^u$, we have
$$
f(B(x',t')) +f(b) \ge  f(B(x',t))+f(B(x,t')).
$$
If $f$ is also increasing, as is $f(u)=e^u$, we obtain
$$
f(B(x',t')) +f(B(x,t)) \ge f(B(x',t')) +f(b) \ge
f(B(x',t))+f(B(x,t')),
$$
which is (\ref{A-sums}).}

\begin{theorem}
\label{thm:main2} If $a: \mathbb{R}^+ \rightarrow \mathbb{R}^+$  is
a uniformly continuous log-concave function then the solution $\tP$ to ${\cal G}(Q)=Q$  is
$\mbox{TP}_2$.
\end{theorem}

\proof Let $\tP_0=0$ and $\tP_{m+1}={\cal G}(\tP_m)$ for $m=1,2,\ldots$. We first claim that the functions $\{\tP_m\}_{m=0}^\infty$ are nonnegative,
increasing in both $x$ and $t$, and $\mbox{TP}_2$. Clearly the
claim is true for $\tP_0$. Assume the claim is true for $\tP_m$ for
some $m \ge 0$. By Lemma \ref{R-is-increasing}  ${\cal
F}(\tP_m)$ and $\tP_{m+1}$ are nonnegative and increasing in $x$
and $t$, and by Lemma \ref{integral-TP2}  ${\cal F}(\tP_m)$ is
$\mbox{TP}_2$, and now Lemma \ref{plus1-TP2} shows that $\tP_{m+1}={\cal F}(\tP_m)+1$ is $\mbox{TP}_2$, thus completing
the proof of the claim. That $\tP$ is $\mbox{TP}_2$ follows, as
$\tP_m$ converges uniformly to $\tP$ by Corollary
\ref{sec:final-corollary}, and the $\mbox{TP}_2$ property is
preserved under (even) pointwise limits. \qed

We may now prove Theorem \ref{main}:

\noindent \textit{Proof of Theorem \ref{main}.} By Lemmas \ref{R-is-bounded} and \ref{lem:F-is-ctype}, when $a$ is a probability, the operator ${\cal G}$ is of $1$-contraction type and $e^t$-bounded. By Corollary \ref{sec:final-corollary} the unique solution $\tP$ to ${\cal G}(Q)=Q$ exists, is continuous and satisfies $0 \le \tP(x,t) \le e^t$.
Now recalling (\ref{tP=etP}), that is, that $\tP(x,t)=e^tP(x,t)$, we have $P(x,t) \in [0,1]$. Theorem \ref{thm:main2} gives that $\tP$ is $\mbox{TP}_2$, hence $P$ is $\mbox{TP}_2$ by Property 5 of Lemma~\ref{lem:tp2}. \qed

The next lemma modifies the result
in \citet{Samuel70}, page 158, to
continuous ammunition and time, thus allowing us to show Theorem
\ref{thm:A}.
\begin{lemma}
\label{esterslemma} Let $R(x,z,t)$ be a positive continuous function of $x
\ge z \ge 0$ and $t \ge 0$ such that the ratio $R(x,z,t)/R(x,z',t)$ is
increasing in $t$ whenever  $x \ge z' \ge z \ge 0$, and let
$$
k(x,t)=\min\{y \in [0,x] : R(x,y,t)=\max_{0 \le z \le x}R(x,z,t)\}.
$$
Then $k(x,t)$ is decreasing in $t$.
\end{lemma}

\proof For the sake of contradiction, assume that there exist $t_1<t_2$ and
$z_1<z_2 \le x$ such that $k(x,t_1)=z_1<z_2=k(x,t_2)$. Then the first
equality gives $R(x,z_1,t_1)/R(x,z_2,t_1) \ge 1$, as the maximum of
$R(x,z,t_1)$ is attained at $z=z_1$, whereas the second equality gives
$R(x,z_1,t_2)/R(x,z_2,t_2)<1$ as $k(x,t_2)$ is the minimum over all $y$
which attain the maximum of $R(x,z,t_2)$. This contradicts the assumption
that $R(x,z,t)/R(x,z',t)$ is increasing in $t$. \qed

\bigskip

\noindent \textit{Proof of Theorem \ref{thm:A}.}  By Theorem \ref{main},
the solution $P(x,t)$ of (\ref{P}) is $\mbox{TP}_2$; also from (\ref{P}), $P(x,t) \ge e^{-t}$ for all $(x,t) \in \mathbb{R}^+ \times \mathbb{R}^+$.
In view of definition
(\ref{TP2-def}), the ratio $P(x',t)/P(x,t)$ is increasing in $t$
whenever $x<x'$, or
$$
\frac{P(x-z,t)}{P(x-z',t)}\,, \quad \mbox{and therefore} \quad
\frac{a(z)P(x-z,t)}{a(z')P(x-z',t)}\,,
$$
is increasing in $t$ when $0\le z<z' \le x$. Applying Lemma \ref{esterslemma}
with the continuous function $R(x,z,t)= a(z)P(x-z,t)$ now yields that
\begin{equation}\label{eq:defK}
K(x,t)=\min\{y \in [0,x]: a(y)P(x-y,t) = \max_{0 \le z \le x}a(z)P(x-z,t)
\}
\end{equation}
is decreasing in $t$, and [A] holds. \qed

We now consider [C] for the Bomber.

\begin{theorem}
\label{C-dcb}
If $a: \mathbb{R}^+ \rightarrow \mathbb{R}^+$  is a strictly log-concave function, then conjecture [C] holds for the doubly continuous Bomber Problem.
\end{theorem}

\proof For $0 \le y < x$ and $\delta=x-y$ note that [C], that is, $x-K(x,t) \ge y-K(y,t)$, is equivalent to
\begin{equation} \label{C-1} K(y+\delta,t) \le K(y,t)+\delta. \end{equation}
Inequality (\ref{C-1})
holds trivially for all $y$ such that $K(y+\delta,t) \le \delta$.
To handle the case $K(y+\delta,t) > \delta$, define
$$
b(x,y)=\frac{a(x)}{a(y)}
\quad \mbox{and} \quad
G_x(y,t) = a(y)P(x-y,t),
$$
the latter expression being the conditional probability of survival upon expending an amount~$y$ of
ammunition at an encounter at time $t$, when the total amount $x$ is
available, and then proceeding optimally thereafter. With these
definitions in place we obtain \begin{equation} \label{Gy+delta}
G_{y+\delta}(v+\delta,t) = a(v+\delta)P(y-v,t)=
\frac{a(v+\delta)}{a(v)}G_y(v,t)=b(v+\delta,v)G_y(v,t). \end{equation} Rewriting
(\ref{Gy+delta}) and using that $K(x,t)$ is optimal, we have \bea G_y(v,t) &=&
\frac{G_{y+\delta}(v+\delta,t)}{b(v+\delta,v)} \nonumber \\
&\le& \frac{G_{y+\delta}(K(y+\delta,t),t)}{b(v+\delta,v)}
\nonumber \\
&=& \frac{G_{y+\delta}([K(y+\delta,t)-\delta]+\delta,t)}{b(v+\delta,v)}
\nonumber \\
\label{C-2} &=&
G_y(K(y+\delta,t)-\delta,t)\frac{b(K(y+\delta,t),K(y+\delta,t)-\delta)}{b(v+\delta,v)},
\ena where in the third expression we use $K(y+\delta,t) > \delta$, and
obtain the last equality by applying (\ref{Gy+delta}) with
$v=K(x,t)-\delta$.

By part \textit{1} of Lemma \ref{lem:tp2} $a(y)$ is strictly $\mbox{TP}_2$, which
is equivalent to $b(x+\delta,x)$ being a strictly decreasing function of $x$ for $x \ge 0$.
Therefore, for $v$ satisfying $v
<K(y+\delta,t)-\delta$, we have
$$
\frac{b(K(y+\delta,t),K(y+\delta,t)-\delta)}{b(v+\delta,v)} <  1.
$$
Since $K(x,t)$ is the minimal amount of ammunition that can be used to
proceed optimally, if the smallest maximizer of the left
hand side of (\ref{C-2}) is attained at a value of $v=K(y,t)$ satisfying $v+\delta
<K(y+\delta,t)$, the ratio of the $b$-functions on the right hand side of
(\ref{C-2}) would be strictly less than 1.  In this case, by (\ref{C-2}),
$K(y+\delta,t)-\delta$  would be a better choice for $v$ than the optimal amount
$K(y,t)$. Note that by definition $K(y+\delta,t) \le y+\delta$, so $K(y+\delta,t)-\delta \le y$, and hence is a feasible quantity to expend, given ammunition $y$. Having reached a contradiction, we conclude
$$
K(y,t)+\delta \ge K(y+\delta,t),
$$
which is (\ref{C-1}). \qed
\ignore{Clearly, $b(x,y)>1$ for $0 \le y < x$.
It seems that one should be able to prove that this inequality is strict.}

\section{The Frail Fighter}
\label{sec:Ffighter}
In this section we consider the Frail Fighter,
whose goal is to shoot down as many of the enemy planes as possible.
Recall that if the Frail Fighter does not shoot down an encountered enemy, then he gets shot
down with probability~1. In this section letting $\otimes$ be given by case $(F_0)$ of (\ref{twostars}),
the dynamic programming equation
(\ref{eq:LPI}) for the Frail Fighter can be written $\tN(x,t)={\cal G}(x,t,\tN)$ for ${\cal G}(x,t,Q)$ given by (\ref{def:calG}) and $\tN(x,t)=e^t N(x,t)$. Analogous to the previous section, in this and the next section $a(y)$ can be interpreted as the probability that the Fighter destroys an enemy using an amount~$y$ of ammunition, although our results hold for more general~$a$.

The following result is parallel to Lemma \ref{integral-TP2}.
\begin{lemma}
\label{le:maxTPpreserves} If $e^{-t}Q(x,t)$ is nonnegative and
increasing in $x$ and $t$ then so is $e^{-t}{\cal G}(Q)$. If $a:
\mathbb{R}^+ \rightarrow \mathbb{R}^+$ is log-concave and $Q$ is
$\mbox{TP}_2(x,t)$ then so is ${\cal G}(Q)$.
\end{lemma}
\proof It is clear that $e^{-t}{\cal G}(Q)$ is nonnegative and increasing in $x$ whenever $e^{-t}Q$ is. Next,
\bea \label{partemtG}
\frac{\partial e^{-t}{\cal G}(Q)}{\partial t}&=&\sup_{0 \le y \le x}a(y)[1 + e^{-t}Q(x-y,t)]-e^{-t}{\cal G}(Q)\\
\nonumber &=&f(t)-\int_0^t f(s)e^{-(t-s)}ds \quad \mbox{where} \quad f(t)=\sup_{0 \le y \le x}a(y)[1 + e^{-t}Q(x-y,t)].
\ena
Since $e^{-t}Q(x,t)$ is increasing in $t$, $f(t) \ge f(s)$ for all $s \in [0,t]$, and
$$
f(t) \ge (1-e^{-t})f(t) = \int_0^t f(t) e^{-(t-s)}ds \ge \int_0^t f(s) e^{-(t-s)}ds.
$$
Hence the partial derivative (\ref{partemtG}) is nonnegative, and $e^{-t}{\cal G}(Q)$ is increasing in $t$.

As $Q(x,t)$ is $\mbox{TP}_2$, so is  $e^{-t}Q(x,t)$ by Lemma \ref{lem:tp2}.
Now, by Lemma \ref{plus1-TP2}, $1+e^{-t}Q(x,t)$ is
$\mbox{TP}_2(x,t)$. As in the proof of Lemma \ref{maxconvpreserves},
letting $p>0$ and $I(x,v)={\bf 1}(0 \le v \le x)$, write
$$
\int_0^x  a(y)^p[1+e^{-s}Q(x-y,s)]^p dy = \int_{-\infty}^\infty  a(x-v)^p
I(x,v) [1+e^{-s}Q(v,s)]^p dv,
$$
which is $\mbox{TP}_2(x,s)$, as in the proof of Lemma \ref{maxconvpreserves}.
Taking $p^{th}$ roots and letting
 $p \rightarrow \infty$ we conclude that $\sup_{0 \le y \le
x}a(y)[1+e^{-s}Q(x-y,s)]$, and therefore also $\sup_{0 \le y \le
x}a(y)[1+e^{-s}Q(x-y,s)]e^s$, are $\mbox{TP}_2(x,s)$.
Now
$$\int_0^t \sup_{0 \le y \le
x}a(y)[1+e^{-s}Q(x-y,s)]e^s ds=\int_{-\infty}^\infty \sup_{0 \le y \le
x}a(y)[1+e^{-s}Q(x-y,s)]I(s,t)e^sds,
$$
and the final claim follows. $\qed$

The proof of Theorem \ref{thm:main3}, which proves the first claim of Theorem \ref{thm:1.3}, is
similar to the proof of Theorem \ref{thm:main2} and is omitted.

\begin{theorem}
\label{thm:main3} If $a: \mathbb{R}^+ \rightarrow \mathbb{R}^+$  is
a uniformly continuous log-concave function then the solution $\tN(x,t)$ to ${\cal G}(Q)=Q$  is
$\mbox{TP}_2$, as is  $N(x,t)=e^{-t}{\overline N}(x,t)$.
\end{theorem}

The Frail Fighter's optimal policy~$K(x,t)$ satisfies
$$
K(x,t)=\min\{y \in [0,x]:  a(y)[1+N(x-y,t)]=  \sup_{0 \le z \le
x}a(z)[1+N(x-z,t)]\}.
$$
Applying Lemma \ref{esterslemma} with $R(x,z,t)=a(z)[1+N(x-z,t)]$ we
obtain [A] for the Frail Fighter, thus proving the last claim of Theorem
\ref{thm:1.3}.

We now prove that [C] holds for the Frail Fighter.
\begin{theorem}
If $a: \mathbb{R}^+ \rightarrow \mathbb{R}^+$  is a strictly log-concave function, then conjecture [C] holds for the Frail Fighter.
\end{theorem}

\proof The argument is nearly identical to the proof of Theorem \ref{C-dcb}, and
again it suffices to only consider $y$ for which $K(y+\delta,t) > \delta$.
Let
$$
b(x,y)=\frac{a(x)}{a(y)}
\quad \mbox{and} \quad
G_x(y,t) = a(y)[1+N(x-y,t)],
$$
the latter expression being the expected `score' for the Fighter expending $y$, conditional
upon an encounter at time $t$ when a total amount $x$ is available, and then
proceeding optimally thereafter. With these
definitions in place we obtain
$$
G_{y+\delta}(v+\delta,t) = a(v+\delta)[1+N(y-v,t)]=
\frac{a(v+\delta)}{a(v)}G_y(v,t)=b(v+\delta,v)G_y(v,t).
$$
The same lines of argument of Theorem \ref{C-dcb} from (\ref{Gy+delta}) on may now
be applied to yield the conclusion.  $\qed$

\section{The Invincible Fighter}
\label{sec:Ifighter} In this section we consider an Invincible
Fighter who cannot be shot down, whose goal is to shoot down as many of the enemy planes as
possible. Hence, in this section we take ${\cal G}$ given by (\ref{def:calG}) and $\otimes$ by case $(F_1)$ of (\ref{twostars}), and the relevant equation may be written $ \tN(x,t)= {\cal G}(x,t,\tN)$. The role of the next lemma, and its proof, are parallel to that of Lemma
\ref{maxconvpreserves} in Section \ref{sec:bomber}.
\begin{lemma}
\label{maxconvpreserves5} If $f(x)$ and $g(x)$ are log-concave functions,
then so is
$$h(x) = \sup_{0 \le y \le x}[f(y)g(x-y)].
$$
If $f(x)$ and $g(x)$ are concave functions, then
so is
$$
h(x) = \sup_{0 \le y \le x}[f(y)+g(x-y)].
$$
\end{lemma}
\proof By Lemma \ref{lem:tp2} part \textit{2}, the convolution of log-concave functions is
log-concave. Writing
$$
h(x)=  \lim_{p \rightarrow \infty} \left(\int_0^x[f(y)g(x-y)]^pdy\right)^{1/p},
$$
and noting that log-concavity is preserved when taking powers and
under pointwise limits, we conclude that $h(x)$ is log-concave, thus
proving the first claim. To prove the second claim, write $h(x) = \log
\sup_{0 \le y \le x}[e^{f(y)}e^{g(x-y)}]$ and apply the first
result to the log-concave functions $e^{f(y)}$ and $e^{g(x)}$. $\qed$

The proof of the following theorem uses the second part of Lemma
\ref{maxconvpreserves5} to show that concavity is preserved under
the operation of $(F_1)$ of (\ref{twostars}), and then an iteration and limit
argument as in the proof of Theorem \ref{thm:main2}. The details are omitted.

\begin{theorem}
\label{thm:main2D} If $a: \mathbb{R}^+ \rightarrow \mathbb{R}^+$ is a uniformly continuous concave function then the solution
$\tN(x,t)$ to the equation ${\cal G}(Q)=Q$, for ${\cal G}$ given in (\ref{def:calG}) for case $(F_1)$ of (\ref{twostars}),
is concave in $x$ for every $t$, as is $N(x,t)=e^{-t}\tN(x,t)$.
\end{theorem}

\ignore{
\proof Let $\tN_m,m=0,1,\ldots$ be given by (\ref{defQm}) with
$\tN_0=0$. Clearly $\tN_0(x,t)=0$ is concave in $x$ for every $t$.
Assuming this property holds for $\tN_m$ for some $m=0,1,\ldots$,
Lemma \ref{maxconvpreserves5} shows that $(a*\tN_m)(x,s)$ is concave
in $x$ for all $s$. Integrating in $s$, we conclude $\tN_{m+1}(x,t)$
is concave in $x$ for every $t$.

As $\tN_m$ converges uniformly to $\tN$ by Corollary
\ref{sec:final-corollary}, and since concavity is preserved under
(even) pointwise limits, the limiting function $\tN$ is concave in
$x$ for every $t$. Clearly the same then must hold for $N$. \qed}

The following lemma is similar to Lemma~\ref{esterslemma}.
\begin{lemma}
\label{ester-type-lemma} If $R(x,y)$ is a continuous function and $R(x,y')-R(x,y)$ is increasing in $x$ for
all $y'>y$, then
$$
k(x)=\inf \{y: R(x,y)=\sup_{0 \le z \le x} R(x,z)\}
$$
is increasing in $x$.
\end{lemma}

\proof To reach a contradiction, assume there exist $x<x'$ such that
$y'=k(x)>k(x')=y$. Note then that by the definition of $k$ we have
$$
R(x,y')=R(x,k(x)) > R(x,k(x'))=R(x,y) \,\,\, \mbox{and} \,\,\, R(x',y')
=R(x',k(x)) \le R(x',k(x'))=R(x',y),
$$
contradicting the assumption on $R$. $\qed$
\bigskip

We may now provide the proof of Theorem \ref{thm:1.4}.

\noindent\textit{Proof of Theorem \ref{thm:1.4}.}  For fixed $t$, let
$$
R(x,y)=a(y)+N(x-y,t).
$$
Note that $R(x,y')-R(x,y)$ being increasing in $x$ for all $y'>y$ is
equivalent to $N(x-y')-N(x-y)$ being increasing in $x$ for all
$y'>y$, which is equivalent to $N$ being concave in $x$. Applying Theorem \ref{thm:main2D}
and Lemma \ref{ester-type-lemma} completes the proof. $\qed$

The optimal policy for the Invincible Fighter is given by
\begin{equation} \label{K-for-IF}
K(x,t)=\arg \max_{0 \le y \le x}[a(y)+N(x-y,t)].
\end{equation}
\begin{proposition}
If the uniformly continuous function $a:\mathbb{R}^+ \rightarrow \mathbb{R}^+$ is strictly concave
then the optimal policy $K(x,t)$ for the Invincible Fighter is unique, that is, the maximum in (\ref{K-for-IF})
 is attained uniquely.
\end{proposition}
\proof Theorem \ref{thm:main2D} gives that $N(x,t)$ is concave in
$x$ for every $t$, and hence $N(x-y,t)$ is concave in $y$ for all
fixed $x$ and $t$. Since $a(y)$ is strictly concave in $y$, so is the
function $a(y)+N(x-y,t)$, and hence its maximum is attained
uniquely. $\qed$

The argument we have used here to show uniqueness no longer works in general in the model $(F_u)$ where the Fighter is vulnerable to
attack, where, following (\ref{eq:gLPI}), the function $R(x,y)$ in the proof of Theorem \ref{thm:1.4} is
replaced by
$$
R(x,y)= a(y)+[a(y)+u(1-a(y))]N(x-y,t).
$$
In this case, $R(x,y')-R(x,y)$ being increasing in $x$ for $y<y'$
is no longer
equivalent to the concavity of $N(x,t)$ in $x$ for fixed $t$ when $u \not =1$.

We now prove that [C] holds for the Invincible Fighter.
\begin{theorem}
If the uniformly continuous function $a: \mathbb{R}^+ \rightarrow \mathbb{R}^+$  is strictly concave, then conjecture [C] holds for the Invincible Fighter.
\end{theorem}

\proof
The argument is similar to the proof of Theorem \ref{C-dcb}. As there,
it suffices to only consider $y$ for which $K(y+\delta,t) > \delta$.
Define
$$
b(x,y)=\exp(a(x)-a(y))
\quad \mbox{and} \quad
G_x(y,t) = \exp(a(y)+N(x-y,t)),
$$
the logarithm of the latter expression being the expected `score' for the Fighter expending $y$, conditional
upon an encounter, at time $t$, when a total amount $x$ is available, and then
proceeding optimally thereafter. With these
definitions in place we obtain
$$
G_{y+\delta}(v+\delta,t) = b(v+\delta,v)G_y(v,t),
$$
and the lines of argument of the proof of Theorem \ref{C-dcb} may now be followed
to show the desired conclusion under the condition that $\exp(a(x))$ is strictly log-concave,
that is, that $a(x)$ is strictly concave. $\qed$

\section{Discussion}
\label{sec:disc}
We have demonstrated that [A] holds for the Bomber and Frail Fighter, [B] for the Invincible Fighter, and [C] for all three, so in no model have we shown that [A] and [B] hold together. In addition, we have shown uniqueness of the optimal policy    only     for the Invincible Fighter. But of all the possible conjectures which can be raised in these models, the one most outstanding is that [B] holds for the Bomber, and naturally one may hope to use the ideas developed in Sections \ref{sec:contraction} and \ref{sec:bomber} to settle this.

\ignore{and in
particular, be able to apply the following variant of Lemma \ref{esterslemma}.
\begin{lemma}
\label{esterslemmavariant}
If $f(x,y)/f(x,y')$ is decreasing in $x$ whenever $0 \le y < y' \le x$ then
$$
k(x)=\inf\{y: f(x,y)=\sup_{0 \le z \le x}f(x,z)\}
$$
is increasing in $x$.
\end{lemma}

\proof For the sake of contradiction assume there exist $x<x'$ such that
$k(x)>k(x')$. Note then that
$$
f(x,k(x)) > f(x,k(x')) \quad \mbox{and} \quad f(x',k(x)) \le
f(x',k(x'))
$$
contradict the fact that $f(x,y)/f(x,y')$ is decreasing when $y=k(x')<k(x)=y'$. \qed

Hence, we would be able to conclude that $K(x,t)$ is
increasing in $x$ if for each fixed $t$
$$
\frac{a(y) \tP(x-y,t)}{a(y') \tP(x-y',t)} \quad \mbox{is
decreasing in $x$ when $0 \le y < y' \le x$.}
$$
Writing out the desired inequality gives, for $0 \le x < x'$,
$$
\frac{a(y) \tP(x-y,t)}{a(y') \tP(x-y',t)} \ge \frac{a(y)
\tP(x'-y,t)}{a(y') \tP(x'-y',t)},
$$
or
$$
\tP(x-y,t) \tP(x'-y',t) \ge \tP(x'-y,t)\tP(x-y',t).
$$
In other words, we would like to show that $\tP(x-y,t)$ is
$\mbox{TP}_2(x,y)$, which is equivalent to $\tP(x,t)$ being log-concave in $x$.}

\begin{figure}[h]
\label{fig:logPm}
\begin{center}\scalebox{.4}{\includegraphics{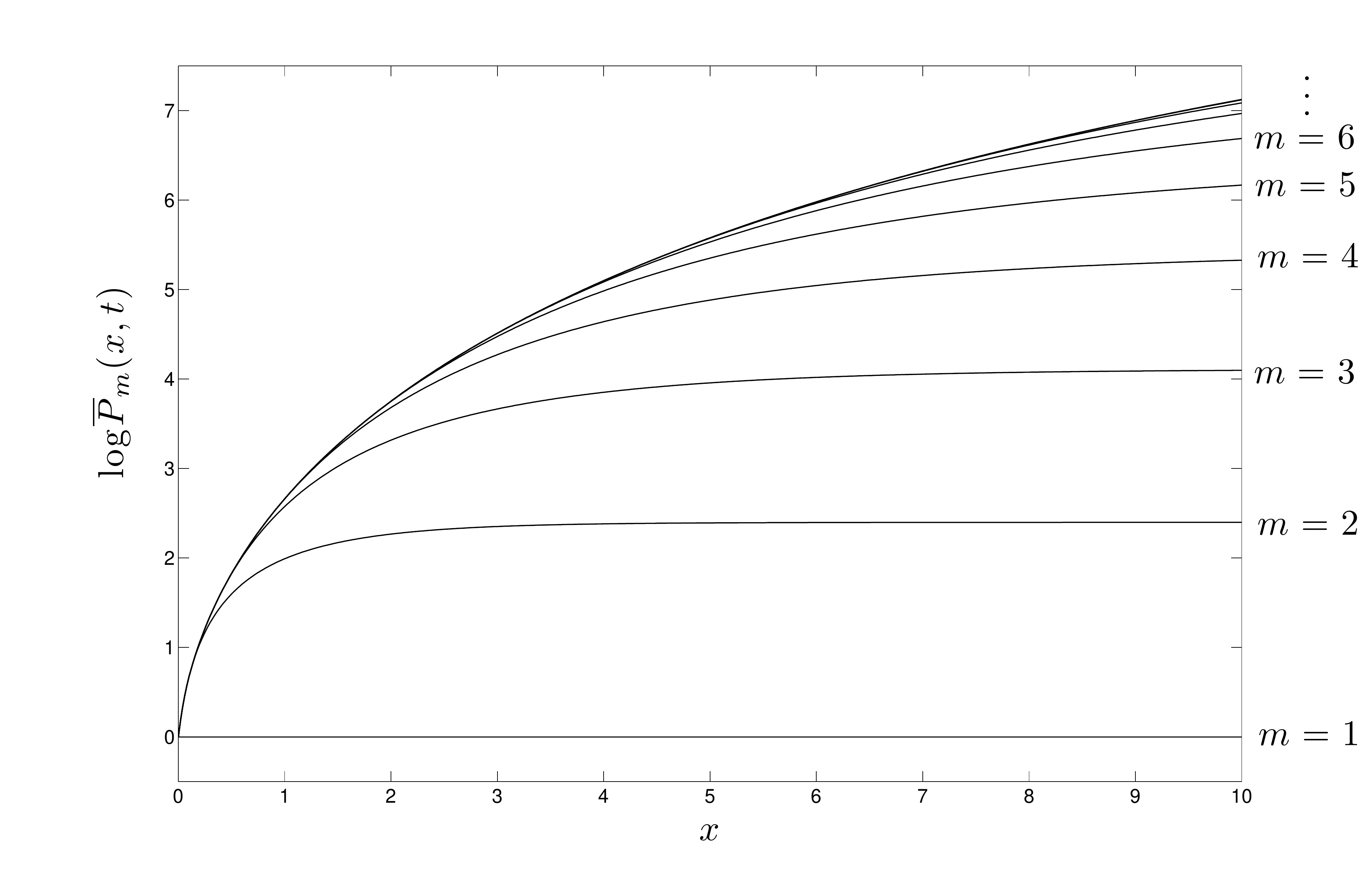}}\\
Figure 1. The functions $\log \overline{P}_m(x,t)$, $m=1,\ldots,10$, for $t=10$.
\end{center}
\end{figure}

\begin{figure}
\label{fig:3dlogPm}
\begin{center}\scalebox{.4}{\includegraphics{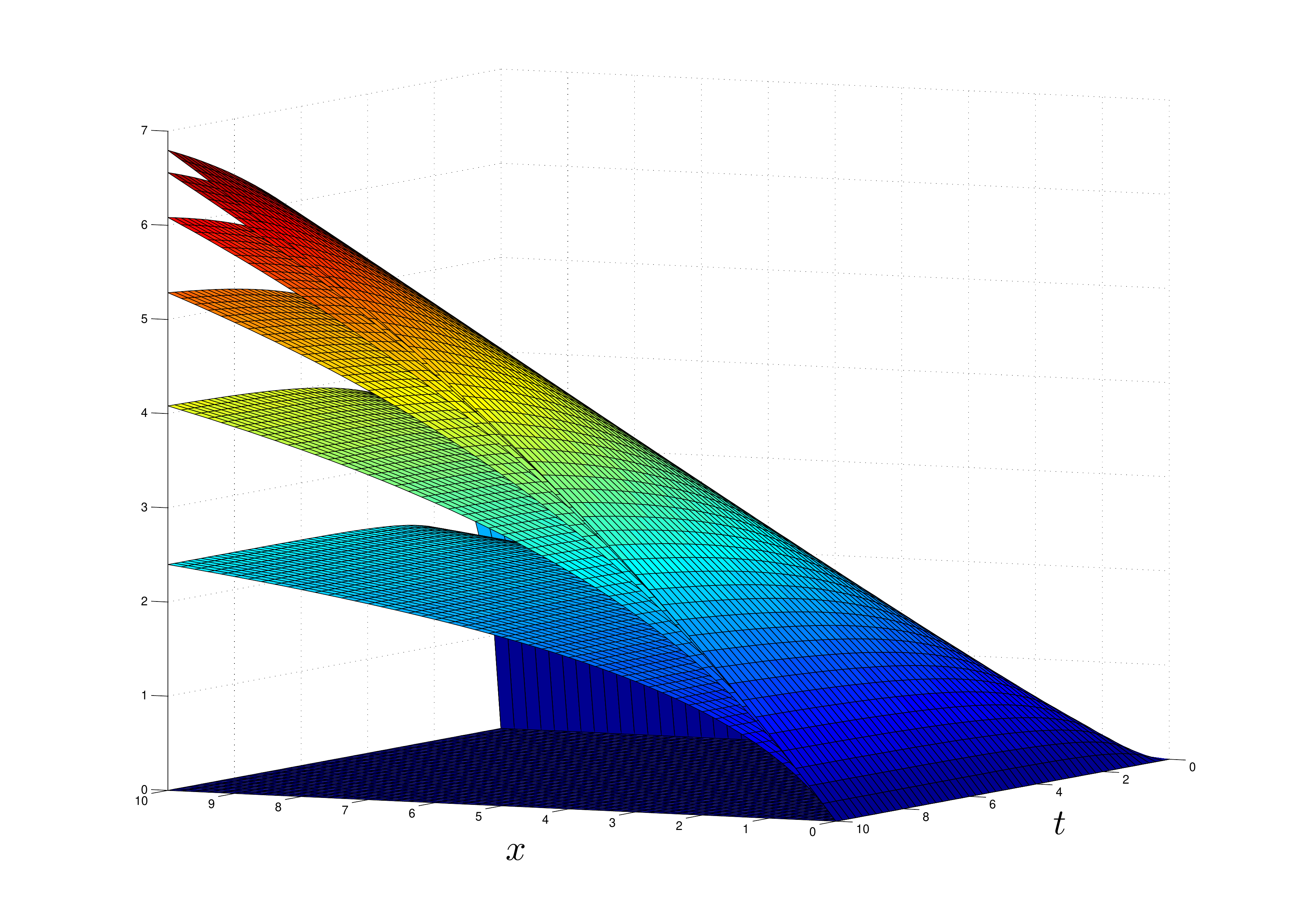}}\\
Figure 2.  The functions $\log \overline{P}_m(x,t)$, $m=1,\ldots, 7$.
\end{center}
\end{figure}

The equivalence of [B] to the log-concavity of $P(x,t)$ in $x$
was noted in \citet{Simons90}. Hence, parallel to the proof of [A], for a
given $Q(x,s)$ consider the function
$$
R(x,t)=\int_0^t (a \otimes Q)(x,s) ds
$$
where $\otimes$ is the max-convolution operator given in case $(Bo)$ of (\ref{twostars}).
By using results in \cite{Prekopa73}, one can show that $R(x,t)$ is log-concave when $a(y)$ and $Q(x,t)$ are log-concave, similarly as to how one shows that the $\mbox{TP}_2$ is preserved.
When $R(x,t)$ is log-concave, for any
nonnegative $x_1,x_2,t_1,t_2$, we have
\begin{equation} \label{Blogcon}
R_{12}^2 \ge R_1 R_2 \quad \mbox{or} \quad R_{12} \ge \left(R_1 R_2 \right)^{1/2},
\end{equation}
where
$$
R_{12}=R(\frac{x_1+x_2}{2},\frac{t_1+t_2}{2}), \quad R_1=R(x_1,t_1) \quad \mbox{and} \quad R_2=R(x_2,t_2).
$$
Now to proceed as in Lemma \ref{plus1-TP2}, for $R(x,t)+1$ to be
log-concave, we require that
$$
\left(R_{12}+1\right)^2 \ge
\left( R_1+1 \right) \left( R_2+1 \right),
$$
and by (\ref{Blogcon}) it suffices that
\begin{equation} \label{Baver}
2 R_{12} \ge
R_1+R_2 \quad \mbox{or} \quad R_{12} \ge \frac{1}{2}(R_1+R_2),
\end{equation}
that is, that $R(x,t)$ is concave. However, unlike the argument in
the proof of Lemma \ref{plus1-TP2}, we cannot conclude (\ref{Baver}) from (\ref{Blogcon}) even though the converse is true. The geometric
mean of $R_1$ and $R_2$ in (\ref{Blogcon}) is in general smaller than the arithmetic mean in (\ref{Baver}). Moreover, numerical evidence
indicates that $R(x,t)$ may fail to be concave.

Nevertheless, as the concavity of $R(x,t)$ is sufficient but not
necessary for the log-concavity of $R(x,t)+1$, log-concavity may yet
be conserved in each step of the iteration, and the limiting
argument pushed through as before. Taking $a(y)$ to be
given in (\ref{a-function}) and considering the first few iterates of
(\ref{recursion}), one has $\tP_0(x,t)=0$ and $\tP_1(x,t)=1$, which
are indeed log-concave. Continuing, as $a(y)$ is increasing we have
$(a \otimes 1)(x,s)=\sup_{0 \le y \le x}a(y)=a(x)$, and therefore
$$
\tP_2(x,t)=\int_0^t (a \otimes \tP_1)(x,s)ds+1=\int_0^t a(x)ds +1 = t a(x)+1.
$$
Since
$$
\frac{\partial^2}{\partial x^2}\log \tP_2(x,t)=\frac{-t(1+t)(1-u)e^{-x}}{(ta(x)+1)^2} \le 0 \quad \mbox{for all $(x,t) \in \mathbb{R}^+ \times \mathbb{R}^+$},
$$
the function $\tP_2(x,t)$ is log-concave as well. Indeed, Figure 1,
which gives plots of $\log \tP_m$ for $m=1,\ldots,10$, $(x,t) \in
[0,10] \times \{10\}$,  and $u=0$, suggests that log-concavity is preserved in these further functions as well. The figure
also indicates that convergence of $\tP_m(x,t)$ for the values
considered is quite rapid, consistent with the exponential bound in
Theorem \ref{thm:Pcmt}. Expanding the range of $t$ back toward the
origin, Figure 2 depicts $\log \tP_m(x,t)$ for $m=1,\ldots,7$ over the
range $(x,t) \in [0,10]^2$, and indicates that these functions are
log-concave in both variables $x$ and $t$.


\def\cprime{$'$}

\end{document}